\documentclass{article}

\usepackage{preamble}

\title{ADVERSARIAL REINFORCEMENT LEARNING: A DUALITY-BASED APPROACH TO SOLVING OPTIMAL CONTROL PROBLEMS}

\author{
        Nan Chen, Mengzhou Liu, Xiaoyan Wang, and Nanyi Zhang \\
		Dept.~of Systems Engineering and Engineering Management\\
        and\\ 
        Centre for Financial Engineering\\
        The Chinese University of Hong Kong, Hong Kong, CHINA
}

\renewcommand{\undertitle}{Accepted by Winter Simulation Conference (WSC) 2025}
\renewcommand{\shorttitle}{Adversarial Reinforcement Learning}

\hypersetup{
pdftitle={A template for the arxiv style},
pdfsubject={q-bio.NC, q-bio.QM},
pdfauthor={David S.~Hippocampus, Elias D.~Striatum},
pdfkeywords={First keyword, Second keyword, More},
}

\begin{document}
\maketitle

\begin{abstract}
We propose an adversarial deep reinforcement learning (ADRL) algorithm for high-dimensional stochastic control problem. Inspired by the information relaxation duality, ADRL reformulates the control problem as a min-max optimization between policies and adversarial penalties, enforcing non-anticipativity while preserving optimality. Numerical experiments demonstrate ADRL’s superior performance to yield tight dual gaps. Our results highlight ADRL's potential as a robust computational framework for high-dimensional stochastic control in simulation-based optimization contexts.
\end{abstract}


\section{INTRODUCTION}
\label{sec:intro}

Stochastic control establishes a powerful framework for modeling and solving decision-making problems in random environments, where uncertainty unfolds and decisions are made sequentially over time. The principle of dynamic programming provides a mathematically elegant characterization of the structure of optimal policies in control problems. However, the practical implementation of this principle is severely constrained by the curse of dimensionality, making it infeasible to solve many high-dimensional applications directly.

Recently an innovative approach emerges in the literature (e.g., \citeNP{he2016} and 
\shortciteNP{hje2018}), leveraging the high expressivity of deep neural networks combined with Monte 
Carlo simulation to address high-dimensional control problems. This deep Monte Carlo 
optimization approach has garnered significant attention for its ability to overcome the 
limitations of traditional methods. For further theoretical developments, diverse applications, 
and numerical experiments, we refer to works such as \shortcite{bcj2019}, 
\shortcite{bgtw2019}, \shortcite{hure2021}. However, some studies (e.g., \citeNP{rs2023}) suggest that 
excessively large neural networks may overfit to future randomness rather than accurately 
estimating it. As a result, the feedback policies trained by this deep Monte Carlo approach 
tend to outperform the original control problem in-sample by implicitly bypassing the essential
adaptiveness restriction. Furthermore, this overfitting issue also causes feedback actions 
constructed using overly large hypothesis spaces to fail to generalize well and perform poorly 
out-of-sample.

The above considerations give rise to an intriguing research problem: how can we assess the optimality of approximate control policies, particularly when these policies are derived using complex approximation schemes such as deep neural networks? In this paper, we propose an adversarial deep reinforcement learning (ADRL) framework for constructing deep neural network-based approximations to the true value functions of multi-dimensional stochastic control problems. This novel approach is rooted in the concept of information relaxation and its associated dual formulation within the stochastic dynamic programming literature (see \citeNP{rogers2007}, \shortciteNP{brown2010}, \citeNP{brown2022}, and \shortciteNP{chen2024}). Drawing inspiration from the principle of strong duality, we find that the dual problem can be reformulated as a zero-sum game between an adversary and an agent. In this game, the adversary seeks to minimize the agent’s expected (penalized) reward by strategically selecting an appropriate dual penalty, while the agent responds by employing non-adaptive policies to maximize its rewards. 
From this perspective, the agent needs to solve deterministic optimization problems in the dual formulation. 

We parameterize the dual penalties using deep neural networks and employ stochastic optimization techniques to develop learning algorithms for identifying the optimal hyperparameters of the dual networks. Both Robbins-Monro (RM) and Kiefer-Wolfowitz (KW) 
type methods are utilized in this framework. Notably, under certain regularity conditions, 
we derive a simple, unbiased gradient estimator for the RM method by leveraging the 
celebrated envelope theorem, which characterizes the differentiability properties of the value function in parameterized optimization problems.

Once the learning process for the dual part is complete, the ADRL algorithm can also generate adaptive control policies. This is achieved by utilizing the learned dual functions as approximations for continuation value functions within the one-step Bellman equation to determine optimal actions. It is important to emphasize that these greedy policies are suboptimal, serving as a lower bound for the true value of the original control problem.
By combining these lower bounds with the upper bounds derived from the dual values, the ADRL algorithm is able to construct effective confidence intervals for the true values of the control problem. Numerical experiments demonstrate that ADRL produces tight confidence intervals even for high-dimensional control problems. These findings underscore the efficiency of the ADRL framework in learning high-quality policies with performance guarantees.

The information relaxation-based duality framework has spurred a rapidly growing body of literature. It has been widely used to demonstrate the near-optimality of certain heuristics by evaluating the tightness of the corresponding dual gap across various applications (see \citeNP{brown2022}, \shortciteNP{chen2024}, and the references therein). However, much of this literature remains silent on how to proceed when the dual gap of a policy under evaluation is found to be loose. Complementing this literature, the main contribution of this paper is the development of a systematic approach for constructing tight dual gaps using deep neural networks, with the ultimate goal of identifying optimal control policies.

\begin{revstart}
Recent advances in adversarial deep learning have primarily focused on efficiently crafting adversarial examples---inputs subtly perturbed to induce model misclassification---to train neural networks and enhance their robustness, ensuring high-confidence performance; see \shortcite{goodfellow2014explaining,huang2017adversarial}. Another prominent technique, Generative Adversarial Networks (GANs) \shortcite{goodfellow2020generative}, employs a min-max formulation to train generator networks through a zero-sum game between generators and discriminators. While our work shares structural similarities with GANs in its use of a min-max optimization framework, it is rooted in a different theoretical foundation and designed for distinct applications. Specifically, our approach leverages information relaxation theory to construct dual values for stochastic control problems. As mentioned earlier, this leads to a systematic approach for constructing tight dual gaps for policy evaluation for solving stochastic control problems. 
\end{revstart}

The remainder of the paper is organized as follows. In Section \ref{sec:control}, we 
review the basics of the SDP duality. Section \ref{sec:main} establishes the ADRL algorithm with discussions on how to implement the Robbins-Monro and Kiefer-Wolfowitz methods. 
Section \ref{sec:numerical} is devoted to one numerical experiment in optimal trade execution.
Section \ref{sec:conclusion} concludes the paper. 
\section{Control Problem and Information Relaxation}
\label{sec:control}

\subsection{Mathematical Formulation}
\label{sec:formulation}

We consider a stochastic optimal control problem in which an agent makes sequential decisions within a random environment over a finite time horizon, denoted as $ \mathbb{T} = \{0, 1, \dots, T\} $. At each time $ t \in \mathbb{T} $, the state of the environment is represented by $\vs_t \in \mathbb{S} \subseteq \mathbb{R}^m $, and the agent selects an action $\va_t \in \mathbb{A}_t \subseteq \mathbb{R}^n $. The agent receives a reward $ r_t(\vs_t, \va_t) $, and the environment transitions to the next state $\vs_{t+1} $ according to the dynamics:
\begin{eqnarray}
\label{dynamic}
\vs_{t+1}=f_t(\vs_{t}, \va_{t}, \xi_{t+1}),
\end{eqnarray}
where $ \xi_{t+1} $ is a random variable taking values in $ \Xi \subseteq \mathbb{R}^d $. The action set $ \mathbb{A}_t $ is characterized by a set of equality and inequality constraints:
\begin{eqnarray}
    \mathbb{A}_t = \{a \in \mathbb{R}^n: g_i(s_t, a) = 0, 1 \le i \le E; \ g_j(s_t, a) \ge 0, 1 \le j \le I\},
\end{eqnarray}
with continuously differentiable constraint functions $\{g_{i}\}_{i=1}^{E}$ and  $\{g_{j}\}_{j=1}^{I}$.

The agent's objective is to maximize the expected total reward over the time horizon $ \mathbb{T} $. Let $ \mathcal{F}_t $ denote the smallest $ \sigma $-algebra generated by the random variables $ {\xi_1, \dots, \xi_t} $. A policy $ \pi = (\pi_0, \dots, 
\pi_{T-1}) $ is said to be \emph{non-anticipative} if $ \pi_t $ is $ \mathcal{F}_t $-measurable, meaning it depends only on the information available up to time $ t $. Let $\Pi$ be the collection of all non-anticipative policies. The agent solves the following optimization problem:
\begin{eqnarray}
\label{value}
V^*(\vs) = \max_{\pi \in \Pi} V^{\pi}(\vs) = \max_{\pi \in \Pi} \E \left[ \sum_{t=0}^{T-1} r_t(\vs_t, \pi_t) + R(\vs_T)  \Big| \vs_0 = \vs \right],
\end{eqnarray}
where $ R(\cdot) $ denotes the terminal reward function. For simplicity, we assume that the random variables $ \{\xi_t\}_{t=1}^T $ are independent. Under this probabilistic structure, the optimal policy $ \pi^* $ is in \emph{feedback} form; that is, $\pi^*_t$ is a function of the current state. In addition, let $ V^*_t(s) $ represent the optimal expected reward achievable from time $ t $ onward, given the current state $ s $. Assume in this paper that 
\begin{revstart}
\begin{asm}
\label{asm:reward}
All the functions $r_{t}(\vs, \va)$, $R(\vs)$, and $f_{t}(\vs, \va, \xi)$ are all bounded and Lipschitz in the sense that there exists a sufficiently large constant $B$ such that, for any $\vs, \vs' \in \mathbb{R}^m$ and $\va, \va' \in \mathbb{R}^n$,
\begin{eqnarray*}
|r_{t}(\vs, \va)-r_{t}(\vs', \va')| \le B(\|\vs-\vs'\|+\|\va-\va'\|),\ 
|R(\vs)-R(\vs')| \le B \|\vs-\vs'\|,
\end{eqnarray*}
and $\sup_{\xi \in \Xi}\|f_{t}(\vs, \va, \xi)-f_{t}(\vs', \va', \xi)\| \le B(\|\vs-\vs'\|+\|\va-\va'\|)$.
\end{asm}
\end{revstart}

\subsection{Information Relaxation and duality}
\label{sec:IRduality}

This paper applies the concept of \emph{information relaxation (IR) based duality} from the literature of stochastic dynamic programming to develop an adversarial reinforcement learning (RL) approach for stochastic optimal control problems. The IR framework relaxes the non-anticipativity constraints, allowing policies to utilize future information, while imposing penalties to account for violations of these constraints. For foundational work on this approach, see \cite{rogers2007}, \shortcite{brown2010}, \cite{brown2022}, and \shortcite{chen2024}.

Let $ \mathbb{A}^T := \mathbb{A}_0 \times \cdots \times \mathbb{A}_{T-1} $ and 
$ \Xi^T := \Xi_1 \times \cdots \times \Xi_T $. We define a function $ z(\va, \xi): \mathbb{A}^T \times \Xi^T \rightarrow \mathbb{R} $ as a \emph{penalty}, where $ \va = (\va_0, \dots, \va_{T-1}) $ is a sequence of actions, and $ \xi = (\xi_1, \dots, \xi_T) $ is a sequence of random perturbations. A penalty function is considered \emph{feasible} if it satisfies $ \mathbb{E}[z(\pi, \xi)] \leq 0 $ for all non-anticipative policies $ \pi \in \Pi $.
Given a feasible penalty function $ z(\cdot, \cdot) $, the corresponding dual value is defined as:
\begin{eqnarray}
\label{dual}
D^z(\vs) := \mathbb{E}\left[\max_{\va = \{\va_t\}_{t=0}^{T-1} \in \mathbb{A}^T} 
\left( \sum_{t=0}^{T-1} r_t(\vs_t, \va_t) + R(\vs_T) - z(\va, \xi) \right)\Big|\vs_0 = \vs\right].
\end{eqnarray}
In this formulation, the optimization is placed inside the expectation, thus relaxing the non-anticipativity constraints. The penalty function $ z $ serves to counterbalance the advantage gained from this relaxation.

The IR literature points out that penalties can be explicitly constructed using a sequence of functions $ W = (W_0, \dots, W_T) $, where $ W_t (\cdot): \mathbb{S} \rightarrow \mathbb{R} $. For $ (\va, \xi) \in \mathbb{A}^T \times \Xi^T $ and $t \in \mathbb{T}$, define 
\begin{eqnarray}
\label{penalty}
z^W_t(\va, \xi) := W_{t+1}(\vs_{t+1}) - \mathbb{E}\left[ W_{t+1}\big(f_t(\vs_t, \va_t, \eta_{t+1})\big)\big|\vs_t \right],
\end{eqnarray}
where $\{\vs_t: t \in \mathbb{T}\} $ is the state trajectory by substituting $\va$ and $\xi$ into the system dynamics (\ref{dynamic}), and $ \eta_t $ shares the same distribution as $ \xi_t $. By summing over all time steps, we can obtain a feasible penalty as follows: 
\begin{eqnarray}
\label{penalty_def}
 z^W(\va, \xi) = \sum_{t=0}^{T-1} z^W_t(\va, \xi).
\end{eqnarray}
Hereafter we refer to $W$ as the \emph{generating function} of IR penalties.  

The strong duality is achievable for finite‑horizon, discrete‑time control problems. Theorem 1 of \cite{rogers2007} and Theorem 2.3 of \shortcite{brown2010} show that
\begin{eqnarray}
\label{strong_dual_spec}
\min_{W}\mathbb{E}\left[\max_{\va=\{\va_t\}_{t=0}^{T} \in \mathbb{A}^T}\left(\sum^{T-1}_{t=0}r_t(\vs_{t}, \va_t) + R(\vs_{T}) -z^W(\va, \xi)\right)\Big|\vs_0 = \vs\right] = V^*(\vs).
\end{eqnarray}
In (\ref{strong_dual_spec}), the minimum on the left-hand side is attained when $W=(V^*_0, \cdots, V^*_T)$, and the optimal policy $\pi^*$ for the primal problem is also optimal for the dual problem. The strong duality relationship (\ref{strong_dual_spec}) presents a systematic approach to assess the optimality of a policy $\pi$. First, we may evaluate the policy by calculating
\begin{eqnarray}
\underline{V}_{t}(\vs)=\mathbb{E}\left[ \sum_{s=t}^{T-1}r_{s}(\vs_s, \pi_{s}) + R_{T}(\vs_{T})\Big| \vs_t=\vs \right]
\end{eqnarray}
for all $0 \le t \le T$. Surely $\underline{V}_{t}(\vs) \le V_{t}(\vs)$ for any $\vs$ because of the sub-optimality of $\pi$. Then, we replace the generic $W$ in (\ref{penalty_def}) by $\underline{V}_{t}$ to construct a penalty $z^{\underline{V}}$ and compute the associated dual value
\begin{eqnarray}
\overline{V}_{0}(\vs)=\mathbb{E}\left[\max_{\va=\{\va_t\}_{t=0}^{T} \in \mathbb{A}^T}\left(\sum^{T-1}_{t=0}r_t(\vs_{t}, \va_t) + R(\vs_{T}) -z^{\underline{V}}(\va, \xi)\right)\Big|\vs_0 = \vs\right]
\end{eqnarray}
From the strong duality (\ref{strong_dual_spec}), we know $V_{0}(\vs) \le \overline{V}_{0}(\vs)$, which implies $0 \le \overline{V}_{0}(\vs)-V_{0}(\vs) \le \overline{V}_{0}(\vs)-\underline{V}_{0}(\vs)$.
When the dual gap  $\overline{V}_{0}-\underline{V}_{0}$ is sufficiently tight, we can conclude that the performance of policy $\pi$ must be very close to the optimality.
This line of work has sparked a rapidly expanding literature that leverages information‑relaxation duality to certify the near‑optimality of a wide range of heuristic policies. More recently, \shortcite{chen2024} develop a supersolution-based approach to improve the optimality for a given policy if we find that its duality 
gap is loose. 
\section{Main results}
\label{sec:main}

\subsection{Adversarial Deep Reinforcement Learning (ADRL)}
\label{sec:adrl}

The min-max formulation of (\ref{strong_dual_spec}) provides a game-theoretic perspective on strong duality. Consider an adversary whose goal is to minimize the agent's expected (penalized) reward by selecting an appropriate $W$ (thereby penalty $z^W$). Meanwhile, the agent seeks to maximize its rewards by employing anticipative policies. From this perspective, the interaction between the agent and the adversary can be framed as a zero-sum game. Based on this observation, we develop an adversarial RL algorithm in this section to learn the solution to (\ref{value}). 

Note that the generating function $W$ of optimal penalties for (\ref{strong_dual_spec}) is typically a highly complex function of state $\vs$, especially in many high-dimensional stochastic control problems. To address this challenge, we leverage the remarkable expressivity of deep neural networks (DNNs) to approximate the adversary's optimal choice for $W$ in this work. In particular, let $\phi_t \in \Phi_t \subset \mathbb{R}^s$ represent the hyperparameters of the network at time $t$, and we parameterize the DNNs as $\{\varrho_{t}(\cdot, \phi_t): t \in \mathbb{T}\}$. One may refer to Appendix \ref{app:dnn} for a brief description on the structure of DNN considered in this paper. Invoking 
(\ref{strong_dual_spec}), we train optimal $\phi_t$ by solving the following minimization problem 
\begin{eqnarray}
\label{opti_phi}
\phi^*=\arg\min_{\phi=\{\phi_t: t \in \mathbb{T}\}}\mathbb{E}\left[\max_{\va=\{\va_t\}_{t=0}^{T} \in \mathbb{A}^T}\left(\sum^{T-1}_{t=0}r_t(\vs_{t}, \va_t) + R(\vs_{T}) -z^{\phi}(\va, \xi)\right)\Big|\vs_0 = \vs\right],
\end{eqnarray}
where $z^\phi$ is a penalty based on the generating function $\{\varrho_{t}(\cdot, \phi_t): t \in \mathbb{T}\}$. 

Structurally, the problem (\ref{opti_phi}) is essentially a stochastic optimization problem. For any hyperparameter $\phi$ in the construction of penalty function, action sequence $\va $, and noise sequence $\xi$, denote 
\begin{eqnarray}
    \label{Y_phi}
    Y(\phi, \va, \xi) := \sum^{T-1}_{t=0}r_t(\vs_{t}, \va_t)+ R(\vs_{T}) -z^{\phi}(\va(\xi), \xi).
\end{eqnarray}    
Here, $\vs=\{\vs_t: t\in \mathbb{T}\}$ is the state sequence determined by the action and randomness pair $(\va, \xi)$. Using this notation, the optimization problem (\ref{opti_phi})
is to maximize the expected value of $\mathbb{E}[\max_{\va}Y(\phi, \va, \xi)]$. 

Inspired by this observation, we develop the following simulation based approach to 
solve it. We generate $n$ sample paths of $\xi$ and collect them into a training dataset $\mathcal{L}^n:=\{\xi^{(1)}, \cdots, \xi^{(n)}\}$, where each $\xi^{(i)}=(\xi^{(i)}_1, \cdots \xi^{(i)}_T)$. Divide $\mathcal{L}^n$ into $m$ disjoint batches: $\mathcal{L}^n = \mathcal{K}_{1}\cup\cdots\cup \mathcal{K}_{m}$. For instance, we break sample paths evenly into $m$ batches: $\mathcal{K}_{1}=\{\xi^{(1)}, \cdots, \xi^{(l)}\}$, $\mathcal{K}_{2}=\{\xi^{(l+1)}, \cdots, \xi^{(2l)}\}$, $\cdots$, $\mathcal{K}_{m}=\{\xi^{((m-1)l+1)}, \cdots, \xi^{(n)}\}$, where $l=n/m$. We use one batch in each iteration of training. Suppose that the first $b-1$ batches $\mathcal{K}_{1}$, $\cdots$, $\mathcal{K}_{b-1}$ have been used in the previous training iterations to yield the current estimate of the hyperparameter $\{\phi^{b-1}_t: t \in \mathbb{T}\}$. The $b$th iteration of ADRL consists of two stages:

-\textbf{Action Stage}: Construct the penalty $ z^{b-1} $ using $\{\varrho_t(\cdot, \phi^{b-1}_t)\}$. The agent solves the inner optimization problem
    \begin{eqnarray}
    \label{action}
    \va^{b}(\xi)&=&(\va^{b}_0(\xi), \cdots, \va^{b}_{T-1}(\xi)) = 
    \arg\max_{\va=\{\va_t\}_{t=0}^{T-1} \in \mathbb{A}^T} Y(\phi^{b-1}, \va, \xi)
    \end{eqnarray}
    along each $\xi \in \mathcal{K}_b$. Therefore, the optimization problem 
    in (\ref{action}) is deterministic. There is a vast research literature about deterministic optimization that we can draw on to solve these problems; see, e.g., \cite{nw2006}. Let us assume for the inner optimization problem, we have 
    \begin{revstart}
    \begin{asm}
    \label{asm:unique}
    (a) The constraint set $\mathbb{A}_t$ is compact and convex for all $t \in \mathbb{T}$. In addition, it satisfies the Slater's condition.\\
    (b) For any given $\phi \in \Phi$, the program (\ref{action}) admits a unique solution for almost every $\xi$.\\
    \end{asm}
    \end{revstart} 

-\textbf{Adversarial Stage}: Update $\phi$ using an estimate $\textrm{grad}^{b-1}_\phi$ of the following gradient:
    \begin{eqnarray}
    \label{eq:gradient}
    \nabla_{\phi}\mathbb{E}\left[\max_{\va=\{\va_t\}_{t=0}^{T-1}}Y(\phi^{b-1}, \va, \xi)\Big|\vs_0 = \vs\right]. 
    \end{eqnarray}
     Update $\phi$ by: $\phi^b = \phi^{b-1} - \gamma_b \cdot \textrm{grad}^{b-1}_\phi$, 
    where $\gamma_b$ is the learning rate. The gradient $\textrm{grad}^{b-1}_\phi$ can be estimated via Robbins-Monro or Kiefer-Wolfowitz methods, detailed in the following subsections.
  
\begin{algorithm}[htbp]
  \caption{Adversarial Deep Reinforcement Learning (ADRL)}
  \label{alg:adrl}
  \begin{algorithmic}[1]
    \Require Training data $\mathcal{L}^n$, network architecture $\{\varrho_t(\cdot, \phi_t)\}$, initial parameters $\{\phi^{0}_t\}$, learning rates $\{\gamma_b\}$
    \For{$b = 1, 2, \dots, m$}
      \State Construct penalty $ z^{b-1} $ using $\{\varrho_t(\cdot, \phi^{b-1}_t)\}$
      \State Select minibatch $\mathcal{K}_b$
      \For{each $\xi \in \mathcal{K}_b$}
        \State Solve (\ref{action}) to find $\va^b(\xi)$
      \EndFor
      \State Compute gradient $\textrm{grad}^{b-1}_\phi$
      \State Update $\phi^b = \phi^{b-1} - \gamma_b \cdot \textrm{grad}^{b-1}_\phi$
    \EndFor
  \end{algorithmic}
\end{algorithm}
We encapsulate the ADRL algorithm in Algorithm \ref{alg:adrl}.  
\begin{revstart}
The ADRL algorithm outputs generating functions $\{\rho_t(\cdot, \phi^m_t), t \in \mathbb{T}\}$. This set of functions can be used to define a greedy policy through the following one-step Bellman equation: for all state $\vs$ at $t \in \mathbb{T}$, we solve
\begin{eqnarray}
\label{policy}
\pi^m_{t}(\vs)&=&\arg\max_{\va_{t} \in \mathbb{A}_{t}}\Big(r_t(\vs, \va_{t})+\mathbb{E}\left[\rho_{t+1}(\vs_{t+1}, \phi^m_{t+1})|\vs_{t}=\vs\right]\Big)\nonumber\\
&=&\arg\max_{\va_{t} \in \mathbb{A}_{t}}\Big(r_t(\vs, \va_{t})+\mathbb{E}\left[\rho_{t+1}(f_{t}(\vs_{t}, a_t, \xi_{t+1}), \phi^m_{t+1})|\vs_{t}=\vs\right]\Big).  
\end{eqnarray}
Note that such calculated $\pi^m$ depends only on the current state $\vs$. Hence, it is an non-anticipative policy. In most practical applications, computing the policy for all possible states in the state space is computationally expensive, and analytical expressions for the expectations are often unavailable. To address this, we calculate $\pi^m$ for a representative subset of states and generalize it to the entire state space using least-squares fitting. In this process, a sample average based on a small-scale simulation of $\xi$ can be used as an approximation for the expectation, facilitating computation; see \cite{bt1996} for further details.\end{revstart} 
When the training data is sufficient, we expect that the generating function 
$\{\rho_t(\cdot, \phi^m_t), t \in \mathbb{T}\}$ learned from ADRL is close to the optimal one; that is, $\rho_t(\cdot, \phi^m_t) \approx V^*_t(\cdot)$. 
By the strong duality relation (\ref{strong_dual_spec}), we know that the greedy 
policy $\pi^m$ defined in (\ref{policy}) is also close to the optimal one. 

\cite{bb2019} develop a general framework in which they derive both a policy and a 
performance bound from the same approximation value function. On the one hand, 
they treat this approximation as the continuation value and apply the greedy method to 
select actions with respect to it. On the other hand, they generate a dual feasible 
penalty, thereby the dual bounds, using the same approximate value. ``Good" 
approximations are thus crucial for this framework to yield tight dual gaps. The ADRL 
algorithm presents a systematic way to construct such ``good" value approximations 
from the dual side of the problem. The penalty $z^m$ constructed from the generating function $\{\rho_t(\cdot, \phi^m_t), t \in \mathbb{T}\}$ should be close to the optimal penalty. Meanwhile, the greedy policy $\pi^m$ should also be nearly optimal, leading to a tight dual gap.

\subsubsection{The Robbins-Monro Method}
\label{sec:RM}
Under some regularity conditions, we can develop unbiased estimators for the gradient (\ref{eq:gradient}). Suppose that the solution to (\ref{action}) exists at $\va^b$ and interchanging expectation and differentiation in (\ref{gradient}) is feasible. Then, 
\begin{eqnarray}
\label{gradient}
&&\nabla_{\phi}\mathbb{E}\left[\max_{\va=\{\va_t\}_{t=0}^{T-1}}Y(\phi^{b-1}, \va, \xi)\Big|\vs_0 = \vs\right]= \mathbb{E}\left[\nabla_{\phi}\max_{\va=\{\va_t\}_{t=0}^{T-1}}Y(\phi^{b-1}, \va(\xi), \xi)\Big|\vs_0 = \vs\right].
\end{eqnarray}
Note that the differentiation operation inside the expectation on the right-hand 
side of the above equality is taken under fixed $\xi$. This pathwise derivative estimation is well-studied in simulation literature; see \cite{g2006}, \cite{ag2007}, and \cite{cf2015} and the reference therein. Hence, 
\begin{eqnarray}
\label{rm_diff}
\nabla_{\phi}\max_{\va=\{\va_t\}_{t=0}^{T-1}}Y(\phi^{b-1}, \va, \xi)=\nabla_\phi  \left(\sum_{t=0}^{T-1} r_t(\vs_t, \va^{b}_t) + R(\vs_T) - z^{b-1}(\va^{b}(\xi), \xi)\right)=\nabla_{\phi}Y(\phi^{b-1}, \va^b(\xi), \xi) 
\end{eqnarray}
provides an unbiased estimator for the gradient on the left-hand side of (\ref{gradient}). 
Fixed $\xi$. Changes in $\phi$ affect $ Y(\phi^{b-1}, \va(\xi), \xi) $ directly through 
$ z^{\phi^{b-1}} $ and indirectly through the optimal solution $\va^{\phi^{b-1}}$ to the inner optimization problem (\ref{action}). However, under certain conditions, the indirect effect is negligible due to the envelope theorem \cite{ms2002}. Using this intuition, we can show that 
\begin{pro}
\label{pro:rm}
Denote 
\begin{eqnarray}
\nabla_{\phi} z^{\phi}_{t}(\va^{\phi}_{t}(\xi), \xi):=\nabla_{\phi}  \varrho_{t+1}(\vs_{t+1}(\va^{\phi}, \xi), \phi_{t+1}) -\mathbb{E}\left[\nabla_{\phi} \varrho_{t+1}\left(f_{t}\left(\vs_{t}(\va^{\phi}, \xi), \va^{\phi}_{t}, \eta_{t+1}\right), \phi_{t+1}\right)\right]
\end{eqnarray}
for $t \in \mathbb{T}$, where $\nabla_{\phi}  \varrho_{t}(\cdot, \phi_t)$ represents the gradient of $\varrho_t$ with respect to $\phi$ and $\{\vs_t(\va, \xi): t \in \mathbb{T}\}$ is obtained by applying action $\va$ and the randomness $\xi$ to the system.
\begin{revstart}
Then, under Assumptions \ref{asm:reward} and \ref{asm:unique}
\end{revstart}
\begin{eqnarray}
\nabla_{\phi}\mathbb{E}\left[\max_{\va=\{\va_t\}_{t=0}^{T-1}}Y(\phi, \va, \xi)\Big|\vs_0 = \vs\right]=\mathbb{E}\left[-\sum_{t=0}^{T-1}\nabla_{\phi} z^{\phi}_{t}(\va^{\phi}_{t}(\xi), \xi)\right]. 
\end{eqnarray}
\end{pro}
\noindent Due to the page limit, we omit the complete proof in this paper and refer readers to the website of the first author (\underline{https://www1.se.cuhk.edu.hk/\textasciitilde nchenweb/Publications.html}) for a proof sketch. According to Proposition \ref{pro:rm}, we construct an unbiased gradient estimator:
\begin{eqnarray}
\textrm{grad}^{b-1}_\phi := -\frac{1}{|\mathcal{K}_b|} \sum_{\xi \in \mathcal{K}_b} \sum_{t=0}^{T-1} \nabla_\phi z^{b-1}_t(\va^b_t(\xi), \xi).
\end{eqnarray}
 This estimator enables the use of the Robbins-Monro method to solve the dual problem. Unlike the Kiefer-Wolfowitz method, this approach does not require solving (\ref{action}) for multiple $\phi$ values.

\subsubsection{The Kiefer-Wolfowitz (KW) Method}

Alternatively, we can use finite difference-based estimates for $\nabla_\phi \mathbb{E}[Y(\phi, \xi)]$. Unlike the Robbins-Monro method, the KW algorithm does not require 
$Y(\phi, \xi)$ to be differentiable. Consider the simultaneous perturbation stochastic approximation (SPSA) algorithm \cite{s1992}. In the $b$th iteration, we generate a random direction $\Delta^b = (\Delta^b_{t,j}: t \in \mathbb{T}, j \in \mathbb{J})$ and perturb $\phi^{b-1}$ by $c_b \Delta^b$. The finite difference approximation is:
\begin{eqnarray}
\label{kw:spsa}
\widehat{\nabla} Y (\phi^{b-1}, \xi) := \frac{\max_{\va}Y(\phi^{b-1} + c_b \Delta^b, \va, \xi) - \max_{\va}Y(\phi^{b-1} - c_b \Delta^b, \va, \xi)}{2c_b} \cdot (\Delta^b)^{-1},
\end{eqnarray}
where $(\Delta^b)^{-1}$ is the element-wise inverse of $\Delta^b$. This estimator requires only two function evaluations, regardless of the dimension of $\phi$. It significantly reducing computational cost compared with the traditional KW estimator \cite{kw1952}. The random directions $\Delta^b$ are i.i.d. with $\mathbb{E}[\Delta^b] = 0$ and finite inverse moments. A common choice is the symmetric Bernoulli distribution ($\pm 1$ with probability $1/2$). The step sizes $\{c_b\}$ and $\{\gamma_b\}$ must satisfy some regularity conditions such as $c_b \to 0$, $\sum_b \gamma_b c_b < \infty$, and 
$\sum_b \gamma^2_b c^{-2}_b < \infty$ to ensure the convergence; see \cite{s1992}, \cite{s2005}, \cite{cf2015}.

\section{Numerical Experiments}
\label{sec:numerical}

Consider an optimal order execution problem in which a trader plans to transact a large 
block of equity over a fixed time horizon with minimum impact costs. It can be viewed as a 
variant of the models proposed in \cite{bert1998}, \cite{ac2001}, and \cite{hw2014}.
Assume that there are $n$ different assets traded in the market, and the trader aims to 
acquire $\bar{\mathbf{R}}$ shares in each of the assets in $T$ periods. The objective of 
the trader is to determine a trading schedule, i.e., how many shares to purchase in each 
period, denoted by $\va_t := \{a_1,  \cdots, a_n\}, \ t = 1,2, \cdots, T$, to minimize the 
associated transaction cost. Let $\mathbf{R}_t \in \mathbb{R}^n$ denote the number of 
remaining orders need to satisfy in each asset at time $t$. Then, a feasible training scheme should satisfy
\begin{eqnarray}
  &&\sum_{t=1}^{T} \va_t = \bar{\mathbf{R}},\ \va_t \geq 0,\ \va_t \in \mathbb{R}^n,\\
  &&\mathbf{R}_{t+1} = \mathbf{R}_{t} - \mathbf{a}_t, \ \mathbf{R}_{1} = 
  \bar{\mathbf{R}}, \ \mathbf{R}_{T+1} = 0,
\end{eqnarray}
where $\va_t \geq 0$ indicate the no-shorting constraint.

To complete the statement of the problem, we specify the price dynamics to be a linear price-impact model, following \cite{he2016}. In particular, assume that price 
$\mathbf{P}_t \in \mathbb{R}^n$ follows
\begin{eqnarray}
  \label{price_dynamic}
  \mathbf{P}_t = \mathbf{P}_{t-1} + \mathbf{A}\mathbf{a}_t + \mathbf{B}\mathbf{X}_t + \epsilon_t, \mbox{ for all } t = 1, \cdots, T
\end{eqnarray}
where $\mathbf{A}\in \mathbb{R}^{n\times n}$ is a positive definite matrix and $\mathbf{B} \in \mathbb{R}^{n \times m}$. Here $\{\epsilon_t, t=1, \cdots, T \}$ is a sequence of white noise with mean zero and covariance matrix $\Sigma_{\epsilon}$. In (\ref{price_dynamic}), the constant matrix $\mathbf{A}$ is used to capture the intensity of the permanent impact: trading the amount of $\mathbf{a}_t$ changes in assets' fundamental values by $\mathbf{A}\mathbf{a}_t$ and this change will last persistently in the future via the iterative relation of $\mathbf{P}$.

In addition, this model allows trader to incorporate predictive signals to extract insights about the stock's future movements, thereby improving the performance of trade execution. The auxiliary process $\mathbf{X}_t \in \mathbb{R}^m$ in (\ref{price_dynamic}) serves this purpose. There are several approaches in the literature for selecting such signals. For example, \cite{bert1998} propose that $\mathbf{X}$ could represent the return of a broader market index, such as the S\&P 500, a factor commonly employed in traditional asset pricing models like the CAPM. Alternatively, $\mathbf{X}$ could be the outputs of an alpha model derived from the trader’s private stock-specific analysis, capturing information that has not yet been fully reflected in market prices. The exact nature of $\mathbf{X}$ is irrelevant in this experiment and thus we abstract its interpretation and assume 
it follows a stationary AR(1) process: $\mathbf{X}_t = \mathbf{C} \mathbf{X}_{t-1} + \eta_t$. The random noise term $\eta_t$ has zero mean and covariance matrix $\Sigma_{\eta}$, independent of $\epsilon_t$. The trader's objective is to minimize
\begin{equation}
  \label{execution:objective}
  \min_{\{\va_{t}, 1 \le t \le T\}} \mathbb{E}\left[\sum_{t=1}^{T}\mathbf{P}^{tr}_{t} \mathbf{a}_{t}\right].
\end{equation}
The state variable of this model can be chosen as $(\mathbf{P}_{t-1}, \mathbf{X}_t, \mathbf{R}_t)$. The control policy is a function of these state variables. 

In general, this problem has no closed form solution; see the discussion in \cite{bert1998} about why this problem becomes intractable due the the non-negative constraints on 
$\{a_t\}$. In our implementation, we set $n = 10$ and $m = 3$, and $T=20$, resulting in a high-dimensional space of control policies: $\mathbb{R}^{23} \rightarrow \mathbb{R}^{10}$. These challenging settings highlight the effectiveness of our ADRL algorithm in 
generating high-quality policies with performance guarantees. For the neural network architecture, we use two hidden layers, each with ReLU as the activation function. We change the number of neurons ($k=10, 50, 100$) in the experiments to test the approximation accuracy of neural networks of different complexity. During training procedure, we apply feature transformation to the state $\mathbf{s_t = (P_{t-1}, X_t, R_t)}$, as this approach has been shown to accelerate convergence and enhance generalization. Specifically, we construct quadratic features as follows: $\left( \mathbf{P_{t-1},\ X_t,\ R_t,\ X_t X_t^{tr},\ R_t R_t^{tr},\ P_{t-1} R_t^{tr},\ X_t R_t^{tr}} \right) \in \mathbb{R}^{263}$.

These transformed features are then fed into the neural network for training. To compute the expectation in penalty construction (\ref{penalty}), we employ a Monte Carlo simulation with a sample size of 2,000. In the Action Stage of the algorithm, we use the Sequential Least Squares Quadratic Programming (SLSQP) method to solve the optimization problem (\ref{action}). In the Adversarial Stage, we apply the Adam optimizer to dynamically adjust the learning rates, $\gamma_b$, utilizing historical gradient information.

\begin{figure}[htb]
\centering
\includegraphics[width=0.4\textwidth]{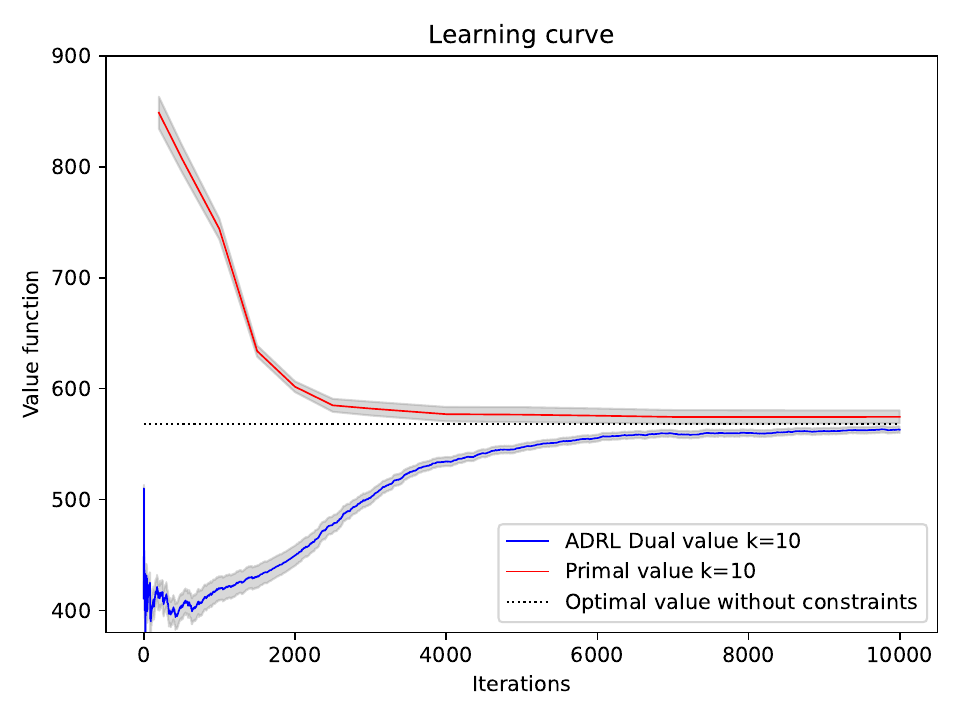}
\caption{ADRL's lower bound (blue) and policy-based upper bound (red) with 95\% confidence intervals (shaded area). All networks trained for $10^4$ iterations on an NVIDIA RTX 4090 GPU (sustained 10\% utilization). The number of neurons $k=10$ in this figure. The total running times is 191 mins.}
\label{fig:adrl_upper_lower_bound}
\end{figure}

The first set of experiments focuses on the case without the no-shorting constraint, $\va \geq 0$. This scenario admits an explicit solution via dynamic programming; see \cite{bert1998} for the analytical expression of the optimal execution strategy and the corresponding optimal cost. For reference, the dotted line in Figure \ref{fig:adrl_upper_lower_bound} represents the value of the optimal cost. The horizontal axis in Figure \ref{fig:adrl_upper_lower_bound} shows the number of iterations. \begin{revstart}
The blue line displays the dual values with the penalty constructed on the generating function $\{\rho_t(\cdot, \phi^m_t), t \in \mathbb{T}\}$ after each iteration.
\end{revstart} The shaded area around the blue line indicates the 95\% confidence interval for these dual estimates. In each iteration, a single sample trajectory of $\xi$ is generated as a mini-batch. 

\begin{revstart}
Theoretically, based on the strong duality relation (\ref{strong_dual_spec}), the obtained dual values should provide lower bounds for the true optimal value. This is supported by the observation that the blue line remains below the dashed line in the figure. Furthermore, Figure \ref{fig:adrl_upper_lower_bound} demonstrates that the dual values generated by our ADRL algorithm converge to the optimal value of the problem (\ref{execution:objective}). As discussed in Section \ref{sec:adrl}, we can construct 
non-anticipative policies using the greedy policy defined on the dual values. Specifically, we take the dual value functions after every 500 iterations and simulate  (\ref{policy}) based on 2,000 samples of ${\xi}$ to evaluate this greedy policy. The red line in Figure \ref{fig:adrl_upper_lower_bound} shows the values associated with these policies. They are all above the dotted line because any policies are suboptimal. However, we can easily see that the greedy policy obtained after 4,000 iterations of ADRL is sufficiently optimal because the policy value is very close to the lower bound provided by dual values. To quantify the policy optimality, define 
$$
\text{gap} := \frac{\text{dual\_value} - \text{primal\_value}}{\text{dual\_value}}.
$$
At 5,000 iterations, the gap reaches 2.82\%. Since the true value lies between the dual and primal values (by strong duality), the relative error of the greedy policy at 5,000 iterations must be less than 2.82\%. As a result, our ADRL algorithm can generate a good policy and meanwhile provide the corresponding dual values as performance guarantee, which presents a systematic way to construct the framework in \cite{bb2019}. \end{revstart} The results for more complex network configurations (k=50 and k=100) exhibit trends similar to those shown in Figure \ref{fig:adrl_upper_lower_bound}. Therefore, we omit the corresponding figures for brevity.



\begin{figure}[htbp]
\centering
\subfloat[ADRL bounds for $k=50$ \label{fig:adrl_con103_k_50}]{
\includegraphics[width=0.3\textwidth]{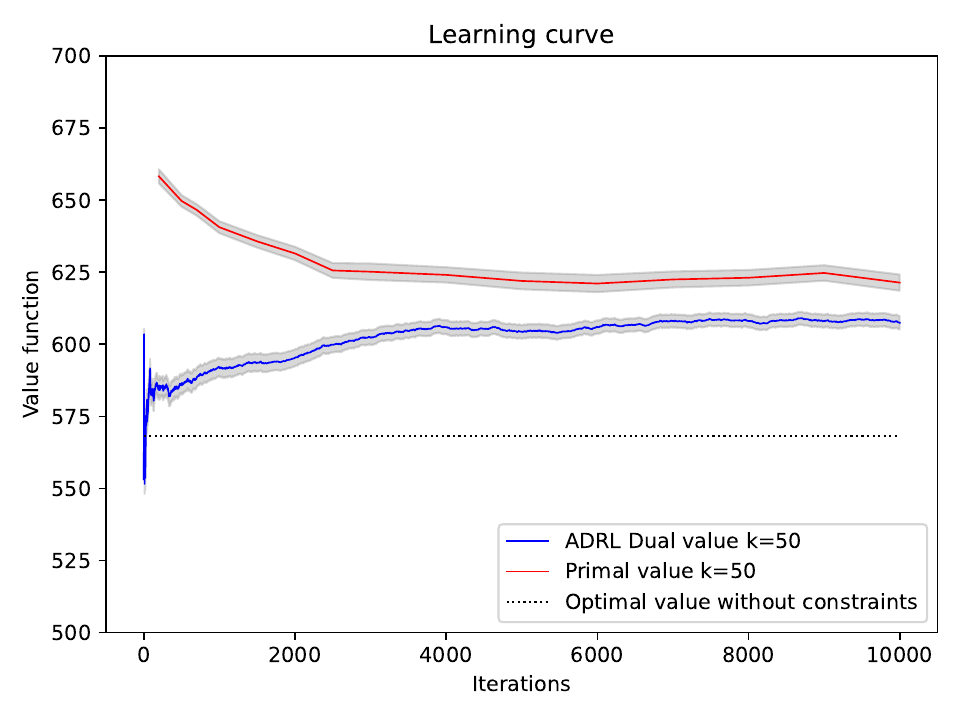}
}
\hspace{0.18\textwidth}
\subfloat[ADRL bounds for $k=100$ \label{fig:adrl_con103_k_100}]{
    \includegraphics[width=0.3\textwidth]{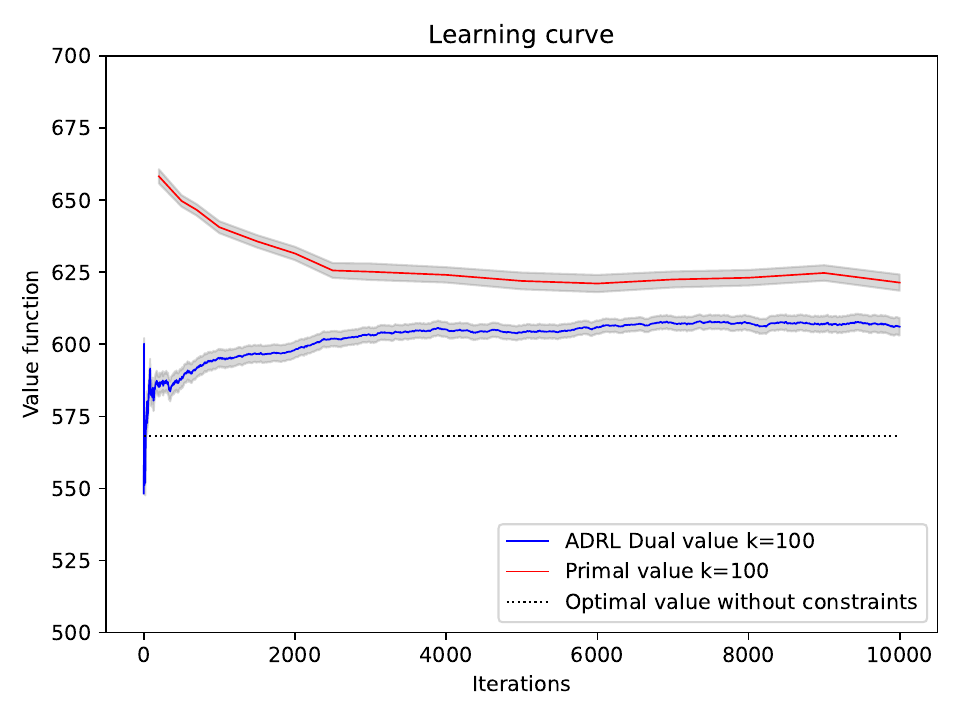}
}
\caption{
ADRL bounds ($k=50, 100$) showing learning curve and policy-based upper bound with 95\% confidence intervals. The dotted line indicates the unconstrained theoretical optimal value. 
}
\label{fig:adrl_combined_results}
\end{figure}

In the second set of experiments, we impose no-shorting constraints on ${\va}$ (i.e., $\va \ge 0$). As noted in \cite{bert1998}, no closed-form solution exists for this case. We continue to examine a ten-asset portfolio setting ($n=10$, $m=3$, $T=20$), consistent with the first experiment. Figures \ref{fig:adrl_con103_k_50} and \ref{fig:adrl_con103_k_100} show that we can still construct effective confidence interval estimates for the problem. After 6,000 iterations, the dual gaps average 1.85\% for $k=50$ and 2.86\% for $k=100$. Based on these tight gaps, we conclude with high confidence that the true value of the control problem in (\ref{execution:objective}) with no-shorting constraints is likely to lie between 607.47 (dual value) and 621.08 (primal value). This value exceeds that of the case where short selling is allowed, which aligns with our expectation that smaller feasible sets lead to higher objective values for minimization problems.

\begin{revstart}
\begin{figure}[htbp]
    \centering
    \includegraphics[width=0.4\linewidth]{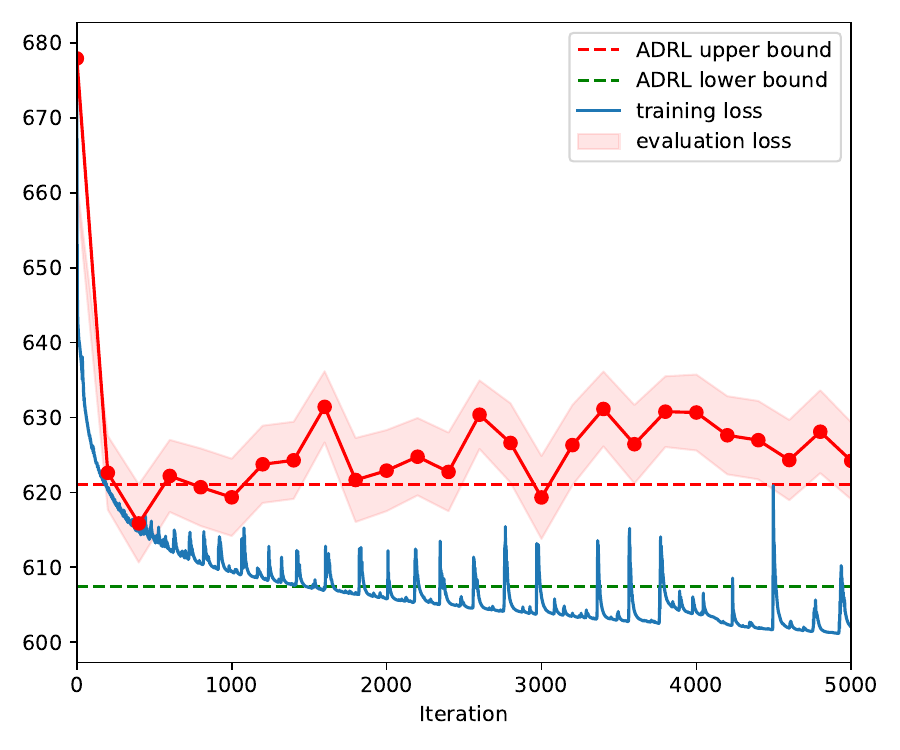}
    \caption{Training and evaluation loss of the ERM-based policy. The blue line shows in-sample training loss; The red dots with shaded bands show the out-of-sample evaluation of the learned policy with 95\% confidence intervals. The two dashed horizontal lines indicate ADRL upper (red) and lower (green) bounds from Figure \ref{fig:adrl_combined_results}.}
    \label{fig:primal}
\end{figure}

Finally, we apply a popular reinforcement learning (RL) algorithm from the literature---the deep empirical risk minimization (DERM) algorithm \cite{he2016}---to the problem with the no-shorting constraint to demonstrate the advantages of working with ADRL. As summarized in \ref{app:derm}, DERM uses in-sample trajectories of the noise process to train policy networks. For this purpose, we simulate 256 samples of $\{(\epsilon_t, \eta_t): 0 \le t \le T\}$ and construct policy networks with two hidden layers, each containing 256 neurons. In each iteration, the samples are used to compute the gradient of the empirical loss on the right-hand side of (\ref{sam_value}) with respect to the hyperparameter $\theta$ for updating the policy networks.

The blue line in Figure~\ref{fig:primal} shows the training loss, which continues to decrease throughout training. However, the evaluation loss (red line) initially improves but eventually worsens, indicating over training, see \cite{rs2023}.
During the initial stages of training, the loss decreases sharply, indicating that the algorithm quickly learns effective policies from the provided data. However, after 2,000 iterations, the algorithm begins to overfit, as evidenced by the training loss dropping below the lower bound provided by ADRL. In other words, the in-sample performance of the trained network surpasses that of the optimal non-anticipative policies. Note that the minimizer $\theta^*$ for (\ref{sam_value}) depends on the entire trajectory ${\vs_t}$. As a result, DERM transitions from learning better non-anticipative policies to effectively predicting future states in the training data. Naturally, such an overfitted policy cannot generalize well to out-of-sample scenarios. As shown by the red dots in the figure, the out-of-sample value of the policies from the latter stage of learning performs worse (higher than the upper bound provided by ADRL).  

A critical problem to resolve this trade off is when to stop the learning. ADRL provides us a systematic approach to address this issue. For instance, when we stop in-sample learning if the training loss lies in between the lower and upper bounds from ADRL, the strong duality relation implies that in-sample performance of the learned policy must be close to the optimality. Figure \ref{fig:primal} illustrates that such stopping rule also leads to good out-of-sample performance. However, we need to acknowledge that, in terms of computational efficiency, the primal method significantly outperforms ADRL, completing 5,000 training iterations in approximately 7 minutes. This advantage is largely attributed to high degree of parallelism. In contrast, ADRL involves repeatedly solving inner optimization problems. How to combine the strength of both algorithms to enhance the scalability of ADRL will be left for future investigation. 
\end{revstart}

\section{Conclusion}
\label{sec:conclusion}

In this study, we develop a novel approach to high-dimensional stochastic control problems by integrating the information relaxation technique with neural network approximations. Our ADRL method simultaneously generates upper and lower bounds for the control problem. Future work will further investigate ADRL's mechanism for mitigating overfitting through adversarial penalties and extend the framework to model-free settings.

\section*{ACKNOWLEDGMENTS}
This research project is partially supported by General Research Fund Scheme (GRF) of Hong Kong Research Grant Council (Project No. 14211023 and 14205824).

\newpage
\bibliographystyle{unsrtnat}
\bibliography{references}  






\newpage
\appendix
\section{APPENDIX} 

\subsection{Neural Networks Architecture}
\label{app:dnn}

In this paper, we employ feedforward neural networks with $N\ge 2$ hidden layers (each containing $k$ neurons) to approximate control policies and value functions. We specify the detailed architecture in this appendix. The network architecture consists of one input layer, one output layer, and $N$ hidden layers. The state vector $\vs$ is encoded into the network through the input layer. More precisely, we can represent it as follows:
\begin{eqnarray}
\label{app:nn}
F(\vs; \phi)=h_{o}\circ\sigma_{N}\circ h_{N} \circ  \cdots  \circ\sigma_{1}\circ  h_{1}(\vs),
\end{eqnarray}
where $h_{1}: \mathbb{R}^m \rightarrow \mathbb{R}^k$, $h_{n}: \mathbb{R}^k \rightarrow \mathbb{R}^k$ ($2 \le n \le N$), and $h_{o}: \mathbb{R}^k \rightarrow \mathbb{R}^d$ are all affine functions. In particular, $h_{1}(\vx)=\mH^1\vx+\vb^1$, $h_{n}(\vx)=\mH^n\vx+\vb^n$, $h_{o}(\vx)=\mH^{o}\vx+\vb^{o}$, where $\mH^1 \in \mathbb{R}^{k \times m}$ (a $k \times m$ matrix), $\mH^n \in \mathbb{R}^{k \times k}$ (a $k \times k$ matrix), $\mH^{o} \in \mathbb{R}^{d\times k}$ (a $d \times k$ matrix), and $\vb^1, \vb^n \in \mathbb{R}^{k}$, for $n=1, \cdots, N$, $\vb^{o} \in \mathbb{R}^d$. $\sigma_{n}: \mathbb{R}^k \rightarrow \mathbb{R}^k$ is a component-wise activation function given by $\sigma_{n}(x_{1}, \cdots, x_{k})=(\sigma(x_{1}), \cdots, \sigma(x_{k}))$ for each $n=1, \cdots, N$.
The hyperparameter $\phi$ in the body text of the paper refers to the matrices $\{\mH^1, \cdots, \mH^N, \mH^{o}\}$ and $\{\vb^1, \cdots, \vb^N, \vb^{o}\}$ used in the construction of network (\ref{app:nn}). 

\subsection{Deep Empirical Risk Minimization Algorithm}
\label{app:derm}
\cite{he2016} developed a deep learning-based approach to directly solves the primal control problem. Its key idea is to approximate the functional dependence of optimal control $\pi^*$ on the state by multilayer feedforward neural networks $\varpi_t(\cdot, \theta_t)$, where $\theta_t$ is the corresponding network hyperparameter. 

To learn the optimal feedback controls, we first repeatedly simulate $n$ i.i.d.
replica of the trajectory of $\xi$: $\mathcal{L}^n:=\{\xi^{(1)}, \cdots, \xi^{(n)}\}$, 
with $\xi^{(i)}=(\xi^{(i)}_1, \cdots, \xi^{(i)}_T)$ for $1\le i \le n$. Then, we may find the optimal hyperparameters through minimizing the following empirical risk
\begin{eqnarray}
\label{sam_value}
v^*(\mathcal{L}^n):=\max_{\{\theta_{t}\}_{t=0}^{T-1}}\frac{1}{n}\sum_{i=1}^{n}\left[\sum^{T-1}_{t=0}r({\vs^{(i)}_{t}}, \varpi_t({\vs^{(i)}_{t}}, \theta_t)) + R({\vs^{(i)}_{T}})\right].
\end{eqnarray}
The standard stochastic gradient descent method with backpropragation can be easily adapted for this purpose. After the in-sample training, we can apply the trained policy network $\varpi(\cdot, \theta^*_t)$ in out of sample to control.

\section*{AUTHOR BIOGRAPHIES}


\noindent {\bf \MakeUppercase{Nan Chen}} is a Professor in the Department of Systems Engineering and Engineering Management at the Chinese University of Hong Kong (CUHK). 
His research interests include financial engineering and FinTech, with a focus on reinforcement learning, quantitative modeling in finance and risk management, Monte Carlo simulation, and applied probability. 
His email address is \email{nchen@se.cuhk.edu.hk} and his website is \url{https://www1.se.cuhk.edu.hk/~nchenweb/}.\\

\noindent {\bf \MakeUppercase{Mengzhou Liu}} is a PhD in the Department of Systems Engineering and Engineering Management at the Chinese University of Hong Kong (CUHK). 
His email address is \email{mzliu@link.cuhk.edu.hk}.\\

\noindent {\bf \MakeUppercase{Xiaoyan Wang}} is a postdoctoral researcher in the Centre for Financial Engineering at the Chinese University of Hong Kong (CUHK). 
Her email address is \email{wangxiaoyan235@link.cuhk.edu.hk}.\\

\noindent {\bf \MakeUppercase{Nanyi Zhang}} is a Master of Science student in Finance in Dept of Financial Mathematics, School of Mathematical Science, Peking University. His email address is \email{2301210046@stu.pku.edu.cn}.\\

\end{document}